\documentclass[11pt]{amsart}

\usepackage{amsmath,amssymb,amscd,amsthm,amsxtra, esint}
\usepackage[dvips]{graphics,epsfig}

\allowdisplaybreaks[2]

\newtheorem{theorem}{Theorem} [section]

\newtheorem{proposition}[theorem]{Proposition}

\newcommand{\noi}{\noindent}
\newcommand{\Z}{\mathbb{Z}}
\newcommand{\R}{\mathbb{R}}

\newcommand{\T}{\mathbb{T}}

\newcommand{\al}{\alpha}
\newcommand{\dl}{\delta}
\newcommand{\Dl}{\Delta}
\newcommand{\eps}{\varepsilon}

\newcommand{\ld}{\lambda}

\newcommand{\s}{\sigma}
\newcommand{\ft}{\widehat}
\newcommand{\wt}{\widetilde}
\newcommand{\cj}{\overline}
\newcommand{\dx}{\partial_x}

\newcommand{\wto}{\rightharpoonup}

\newcommand{\jb}[1]
{\langle #1 \rangle}

\numberwithin{equation}{section}
\numberwithin{theorem}{section}

\begin{document}

%\baselineskip = 20pt
%\date{\today} 

\title
[Cubic NLS below $L^2$]
{On the one-dimensional cubic nonlinear Schr\"odinger equation below $L^2$}

\author{Tadahiro Oh and Catherine Sulem}

\address{
%Tadahiro Oh\\
Department of Mathematics\\
University of Toronto\\
40 St. George St, 
%Rm 6290,
Toronto, ON M5S 2E4, Canada}

\email{oh@math.toronto.edu}

%\address{Catherine Sulem\\
%Department of Mathematics\\
%University of Toronto\\
%40 St. George St, Rm 6290,
%Toronto, ON M5S 2E4, Canada}
 \thanks{C.S. is partially supported by 
            N.S.E.R.C. Grant 0046179-05.}
\email{sulem@math.toronto.edu}



\begin{abstract}
In this paper, we review  several recent results concerning well-posedness of
the one-dimensional, cubic Nonlinear Schr\"odinger equation  (NLS) on the real
line $\R$ and on the circle $\T$ for solutions below the  $L^2$-threshold. 
We point out common results for  NLS on $\R$ and the so-called 
{\it{Wick ordered NLS} } (WNLS) on $\T$, suggesting that WNLS may be an
appropriate model
for the study of solutions below $L^2(\T)$. In particular, in contrast
with  a recent
result of Molinet \cite{MOLI} who proved that the solution map for
 the periodic cubic NLS equation is not weakly continuous from
 $L^2(\T)$ to the space of distributions, we show that this is not
the case for WNLS.  
\end{abstract}

\maketitle

%\tableofcontents

\section{Introduction}

In this paper, we consider the one-dimensional cubic nonlinear Schr\"odinger equation (NLS):
\begin{equation} \label{NLS1}
\begin{cases}
i u_t - u_{xx} \pm |u|^2 u = 0\\
u\big|_{t = 0} = u_0,
\end{cases}
\quad (x, t) \in \T\times \R \text{ or } \R \times \R,
\end{equation}

\noi
where $u$ is a complex-valued function and $\T =\R/2\pi\Z$.
\eqref{NLS1} arises in various physical settings for the description
of wave propagation in nonlinear optics, fluids and plasmas 
(see  \cite{SULEM} for a general review.)
It also arises in quantum field theory as a mean field equation for many body
boson systems.
It is known to be one of the simplest partial differential equations (PDEs) with complete integrability
 \cite{A1, A2, GKP}.

As a complete integrable PDE, \eqref{NLS1}
enjoys infinitely many conservation laws,
starting with conservation of mass, momentum, and Hamiltonian:
\begin{equation} \label{conserved}
N(u) = \int |u|^2 dx, \quad 
P(u) = \text{Im} \int \cj{u} u_x dx, \quad
H(u) = \frac{1}{2}\int |u_x|^2 dx \pm \frac{1}{4} \int |u|^4 dx.
\end{equation}

\noi
In the focusing case (with the $-$ sign), 
\eqref{NLS1} admits soliton and multi-soliton solutions.
Moreover, 
\eqref{NLS1} is globally well-posed in $L^2$
thanks to the conservation of the $L^2$-norm
(Tsutsumi \cite{Tsutsumi} on $\R$ and 
Bourgain \cite{BO1} on $\T$.)

It is also well-known that \eqref{NLS1} is invariant under several symmetries.
In the following, we concentrate on the dilation symmetry and the Galilean symmetry.
The dilation symmetry states that if $u(x, t)$ is a solution to \eqref{NLS1}
on $\R$ with initial condition $u_0$, then
$u^\ld(x, t) = \ld^{-1} u (\ld^{-1}x, \ld^{-2}t)$
is also a solution to \eqref{NLS1} with the $\ld$-scaled initial condition 
$u_0^\ld(x) = \ld^{-1} u_0 (\ld^{-1}x)$.
Associated to the dilation symmetry, 
there is a  scaling-critical Sobolev index $s_c$
such that the homogeneous $\dot{H}^{s_c}$-norm is invariant
under the dilation symmetry. 
In the case of the one-dimensional cubic NLS, 
 the scaling-critical Sobolev index is $s_c = -\frac{1}{2}$.
It is commonly conjectured that a PDE is ill-posed in $H^s$ for $s < s_c$.
Indeed, on the real line, Christ-Colliander-Tao \cite{CCT} showed
that the data-to-solution map is unbounded from $H^s(\R)$ to $H^s(\R)$
for $s<-\frac{1}{2}$. 
The Galilean invariance states that if $u(x, t)$ is a solution to \eqref{NLS1}
on $\R$ with initial condition $u_0$, then
$u^\beta(x, t) = e^{i\frac{\beta}{2}x}e^{i\frac{\beta^2}{4}t}  u (x+\beta t, t)$
is also a solution to \eqref{NLS1} with the initial condition 
$u_0^\beta(x) = e^{i\frac{\beta}{2}x}  u_0 (x)$.
Note that the $L^2$-norm is invariant under the Galilean symmetry.\footnote{The 
Galilean symmetry does not preserve the momentum.
Indeed, $P(u^\beta) = \frac{\beta}{2} N(u) + P(u)$.}
It turned out that this symmetry also leads to a kind of ill-posedness in the
sense that 
 the solution map cannot be {\it smooth} in $H^s$ 
for $s < s_c^\infty = 0$.
Indeed, a simple application of Bourgain's idea in \cite{BO3}
shows that the solution map of \eqref{NLS1}
cannot be $C^3$ in $H^s$ for $s < s_c^\infty = 0$.
See Section 2 for more results in this direction.

Recently, Molinet \cite{MOLI} showed that 
the solution map for  \eqref{NLS1} on $\T$
cannot be continuous in $H^s(\T)$ for $s < 0$.
(See Christ-Colliander-Tao \cite{CCT2} and Carles-Dumas-Sparber \cite{CDS} for related results.)
His argument is based on showing that
the solution map\footnote{Strictly speaking, Molinet's result applies to the flow map,
i.e. for each nonzero $u_0 \in L^2(\T)$, the map: $u_0 \to u(t)$ is not continuous.}  
is not continuous
from $L^2(\T)$ endowed with weak topology
to the space of distributions $(C^\infty(\T))^*$.
Several remarks are in order.
First,  on the real line, 
there is no corresponding result (i.e. failure of continuity of the 
solution map for $s< 0.$) Also,  
the discontinuity in \cite{MOLI} is precisely caused by
$2 \mu(u) u $, where $\mu(u) := \fint |u|^2 dx = \frac{1}{2\pi} \int_0^{2\pi} |u|^2 dx$.

Our main goal in this paper is to  propose  an alternative formulation 
of the periodic cubic NLS below $L^2(\T)$
to avoid this non-desirable behavior.
In particular, we show that this model has properties similar 
 to those of \eqref{NLS1} on the real line
even below $L^2$.
We consider the {\it Wick ordered cubic NLS} (WNLS):
\begin{equation} \label{WNLS1}
\begin{cases}
i u_t - u_{xx} \pm (|u|^2 - 2 \fint |u|^2) u = 0\\
u\big|_{t = 0} = u_0
\end{cases}
\end{equation}

\noi
for $(x, t) \in \T \times \R$.
Clearly,  \eqref{NLS1} and \eqref{WNLS1}
are equivalent for $u_0 \in L^2(\T)$.
If $u$ satisfies \eqref{NLS1} with $u_0\in L^2(\T)$,
then
$v(t) = e^{\mp 2 i \mu(u_0) t} u(t)$ satisfies \eqref{WNLS1}.
However, for $u_0 \notin L^2(\T)$, 
we cannot freely convert solutions of \eqref{WNLS1} into solutions of \eqref{NLS1}.
The effect of this modification can be seen more clearly
on the Fourier side.
By writing the cubic nonlinearity as
$\ft{|u|^2u}(n) = \sum_{n = n_1 - n_2 + n_3}
\ft{u}(n_1)\cj{\ft{u}}(n_2)\ft{u}(n_3)$,
we see that the additional term in \eqref{WNLS1} precisely removes resonant interactions
caused by $n_2 = n_1$ or $n_3$.
See Section 4.
Such a modification does not seem to have a significant effect on $\R$,
since $\xi_2 = \xi_1$ or $\xi_3$
is a set of measure zero in the hyperplane $\xi = \xi_1-\xi_2+\xi_3$
(for fixed $\xi$.)

It turns out that \eqref{WNLS1} on $\T$ shares many common features with
\eqref{NLS1} on $\R$ (see Section 2.)
Equation  \eqref{WNLS1} (in the defocusing case on $\T^2$)
first appeared in the work of Bourgain \cite{BO6, BO7}, 
in the study of  the invariance of the Gibbs measure, 
as an equivalent formulation
of the Wick ordered Hamiltonian equation,
related to renormalization in the Euclidean $\varphi^4_2$ quantum field theory
(see Section 3.) 

There are several results on \eqref{WNLS1}.
Using a power series method, Christ \cite{Christ} proved the local-in-time
existence of  solutions in $\mathcal{F}L^p(\T)$
for $p< \infty$,
where the Fourier-Lebesgue space $\mathcal{F}L^p(\T)$
is defined by the norm $\|f\|_{\mathcal{F}L^p(\T)} = \|\ft{f}(n)\|_{l^p_n(\Z)}$.
In the periodic case, we have $\mathcal{F}L^p(\T) \supsetneq L^2(\T)$ for $p>2$.
Gr\"unrock-Herr \cite{GH} established the same result (with uniqueness) 
via the fixed point argument.

On the one hand,  Molinet's ill-posedness result does not apply to \eqref{WNLS1}
since we removed the part responsible for the discontinuity.
On the other hand, by a slight modification of the argument in Burq-G\'erard-Tzvetkov \cite{BGT},
we see that the solution map for \eqref{WNLS1} is not uniformly continuous below $L^2(\T)$,
see \cite{CO1}.
This, in particular, implies that one cannot expect 
 well-posedness of \eqref{WNLS1}
in $H^s(\T)$ for $s<0$ via the standard fixed point argument.

There are however positive results for \eqref{WNLS1} in $H^s(\T)$ for $s<0$.
Christ-Holmer-Tataru \cite{CHT} established 
 an a priori bound on the growth of (smooth) solutions
in the $H^s$-topology for $s\geq -\frac{1}{6}$.
In Section 4, we show that the solution map for \eqref{WNLS1} 
is continuous in $L^2(\T)$ endowed with weak topology.
These results have counterparts for \eqref{NLS1} on $\R$.

In \cite{CO1}, Colliander-Oh considered the well-posedness question
 of \eqref{WNLS1}
below $L^2(\T)$ with randomized  initial data of the form
\begin{equation}\label{IV}
u_0 (x;\omega) = \sum_{n\in\Z} \frac{g_n(\omega)}{\sqrt{1+|n|^{2\al}}}e^{inx}, 
\end{equation}

\noi
where $\{g_n\}_{n\in\Z}$ is a family of independent standard 
complex-valued Gaussian random variables.
It is known \cite{Z} that
 $u_0 (\omega) \in H^{\al-\frac{1}{2}-\eps}\setminus H^{\al-\frac{1}{2}}$ 
almost surely in $\omega$ for any $\eps > 0$
and that $u_0$ of the form \eqref{IV} is a typical element
in the support of the Gaussian measure
\begin{equation} \label{Gaussian1}
d \rho_\al = Z^{-1}_\al \exp \Big( -\frac{1}{2} \int |u|^2 -\frac{1}{2} |D^\al u|^2 dx\Big) 
\prod_{x \in\T} d u(x),
\end{equation}

\noi
where $D = \sqrt{-\dx^2}$.
In \cite{CO1},   local-in-time solutions were constructed for \eqref{WNLS1}
almost surely (with respect to $\rho_\al$)
in $H^s(\T)$ for each $s > -\frac{1}{3}$
($s = \al - \frac{1}{2}-\eps$ for small $\eps > 0$),
and global-in-time solutions
almost surely 
in $H^s(\T)$ for all $s > -\frac{1}{12}$.
The argument is based on the fixed point argument around the linear solution,
exploiting nonlinear smoothing under randomization on initial data.

The same technique can be applied to study the well-posedness issue of \eqref{WNLS1}
with initial data of the form
\begin{equation}\label{IV2}
u_0 (x;\omega) = v_0(x) + \sum_{n\in\Z} \frac{g_n(\omega)}{\sqrt{1+|n|^{2\al}}}e^{inx}, 
\end{equation}

\noi
where $v_0$ is in $L^2(\T)$.
The initial data of the form \eqref{IV2} may be of physical importance
since smooth data may be perturbed by a rough random noise.
i.e. initial data, which are smooth in an ideal situation,
may be of low regularity in practice due to a noise. 
This is one of the reasons 
that we are interested in having a formulation of NLS below $L^2$.

Another physically relevant issue is the study of \eqref{WNLS1}
with initial data of the form \eqref{IV} when $\al = 0$.
The Gaussian measure $\rho_\al$ then  corresponds
to the white noise on $\T$
(up to a multiplicative constant.)
It is conjectured \cite{Z}
that the white noise is invariant under the flow of the cubic NLS \eqref{NLS1}.
In \cite{OQV}, Oh-Quastel-Valk\'o 
proved that the white noise is a weak limit of probability measures
that are invariant under the flow of \eqref{NLS1} and \eqref{WNLS1}.
Note that the white noise $\rho_0$ is supported 
on $H^{-\frac{1}{2}-\eps}(\T)\setminus H^{-\frac{1}{2}}(\T)$ for $\eps > 0$
(more precisely, on $B^{-\frac{1}{2}}_{2, \infty}$.)
Such a low regularity seems to be out of reach at this point.
Hence, the result in \cite{OQV} implies only a version of ``formal'' invariance of the white noise
due to lack of well-defined flow of NLS on the support of the white noise.
In view of Molinet's ill-posedness below $L^2(\T)$, 
we need to pursue this issue with \eqref{WNLS1} in place of \eqref{NLS1}.
In this respect, the result in \cite{CO1} can be regarded as 
a partial progress toward this goal.

Note that the white noise (i.e. $u_0$ in \eqref{IV} with $\al = 0$ up to multiplicative constant) 
can be regarded as
a Gaussian randomization (on the Fourier coefficients) of the delta function $\dl(x) = \sum_n e^{inx}$.
It is known \cite{KPV5} that 
in considering the Cauchy problem \eqref{NLS1} on $\R$ with the delta function as initial condition,
we have either non-existence or non-uniqueness in $C ([-T, T]; \mathcal{S}'(\R))$.
Moreover, on $\T$, Christ \cite{C2} proved 
a non-uniqueness result of \eqref{WNLS1} in 
the class $C ([-T, T]; H^s(\T))$ for $s< 0$. 
Christ's result states that 
one can not have unconditional uniqueness\footnote{
We say that a solution $u$ is unconditionally unique
if it is unique in $C([0, T];H^s)$
{\it without} intersecting with any auxiliary function space.
Unconditional uniqueness is a concept of uniqueness which does not depend
on how solutions are constructed.
See Kato \cite{KATO}.}
in $H^s(\T)$, $s<0$.
However, this is not an issue 
since, in  discussing well-posedness, 
we usually construct a unique solution in $C ([-T, T]; H^s) \cap X_T$, where $X_T$ is an auxiliary function space
(such as Strichartz spaces or $X^{s, b}$ spaces.)

Lastly, another physical motivation 
for the study of NLS in the low regularity setting
is the localized induction approximation model 
for the flow of a vortex filament.
The filament at time $t$ is given by 
a curve $X(x, t)$ in $\R^3$,
satisfying
\begin{equation} \label{filament}
X_t = X_x \times X_{xx},
\end{equation}

\noi 
where $x$ is the arclength.
Then, under the Hasimoto transform \cite{HA}:
\begin{equation} \label{HASIMOTO}
 u(x, t) = c(x, t) \exp \bigg( i \int^x \tau (y, t) dy \bigg), 
\end{equation}

\noi
where $c(x, t)$ and $\tau(x, t)$ are the curvature and  the torsion of $X(x, t)$,
the transformed function $u$ satisfies the focusing cubic \eqref{NLS1} on $\R$.
Guti\'errez-Rivas-Vega \cite{GRV}
showed that 
a smooth filament can develop a sharp corner in finite time,
which corresponds, under  \eqref{HASIMOTO}, 
to a Dirac delta singularity for $u$ in \eqref{NLS1}.
This necessitates 
the study of NLS in the low regularity setting.

\medskip
This paper is organized as follows.
In Section 2, we compare the results for NLS \eqref{NLS1} on $\R$
and Wick ordered NLS \eqref{WNLS1} on $\T$.
In Section 3, we recall basic aspects of the Wick ordering
and the derivation of \eqref{WNLS1} on $\T^2$ following \cite{BO7}.
In Section 4, we present the proof of the weak continuity of
the solution map for \eqref{WNLS1} in $L^2(\T)$.

\section{NLS on $\R$ and Wick ordered NLS on $\T$}

In this section, we present several results
that are common to  \eqref{NLS1} on $\R$ and 
\eqref{WNLS1} on $\T$.
We show a  summary of these results in Table \ref{TAB1} below. 
This analogy suggests that  
Wick ordered NLS \eqref{WNLS1} on $\T$ is an appropriate  model
to study  when interested in solutions below $L^2(\T)$.

\subsection{Well-posedness in $L^2$}
On the real line, Tsutsumi \cite{Tsutsumi} proved global well-posedness of \eqref{NLS1} in $L^2(\R)$.
His argument is based on the smoothing properties of the linear Schr\"odinger operator
expressed by the Strichartz estimates and the conservation
 of the $L^2$-norm. %It also applied to higher dimensions. 
For the problem on the circle, Bourgain \cite{BO1} introduced 
the $X^{s, b}$ space
and proved global well-posedness of \eqref{NLS1} in $L^2(\T)$.
His argument is based on the periodic $L^4$ Strichartz
and the conservation of the $L^2$-norm.
The same argument can be applied to establish
global well-posedness of \eqref{WNLS1} in $L^2(\T)$.

\subsection{ Ill-posedness in $H^s$ for $s<0$:}

An application of Bourgain's argument in \cite{BO3}
shows that the solution maps for \eqref{NLS1} on $\R$ and \eqref{WNLS1} on $\T$ are not $C^3$ in $H^s$ for $s<0$.
The method consists of examining the differentiability at $\dl = 0$ of the solution map
with initial condition $u_0 = \dl \phi$ for some suitable $\phi$
i.e. differentiability at the zero solution in a certain direction.

On $\R$, Kenig-Ponce-Vega \cite{KPV5} proved the failure of uniform continuity 
of the solution map for \eqref{NLS1} in $H^s(\R)$ for $s<0$ in the focusing case,
by constructing a family of smooth soliton solutions.
In the defocusing case, Christ-Colliander-Tao \cite{CCT}
established the same result
by constructing a family of smooth approximate solutions.
On $\T$, Burq-G\'erard-Tzvetkov \cite{BGT} (also see \cite{CCT}) constructed
a family of explicit solutions supported on a single mode
and showed the corresponding result for \eqref{NLS1}.
By a slight modification of their argument, 
we can also establish the same result for \eqref{WNLS1}. 
It is worthwhile to note that the momentum diverges to $\infty$
in these examples.

The above ill-posedness results 
state that the solution map is not smooth or uniformly continuous
in $H^s$ below $s < s_c^\infty = 0$.
This does not say that \eqref{NLS1} on $\R$ and \eqref{WNLS1} on $\T$
are ill-posed below $L^2$,
i.e. it is still possible to construct continuous flow below $L^2$.
These results instead state that the fixed point argument
cannot be used to show well-posedness of \eqref{NLS1} on $\R$
 and \eqref{WNLS1} on $\T$ below $L^2$,
since such a method would make  solution maps smooth.
Compare the above results with
the ill-posedness result by Molinet \cite{MOLI}
- the discontinuity of the solution map below $L^2(\T)$
for the periodic NLS \eqref{NLS1}.

\subsection{Well-posedness in $\mathcal{F}L^p$:}
Define the Fourier Lebesgue space $\mathcal{F}L^{s, p}(\R)$
equipped with the norm
$\|f\|_{\mathcal{F}L^{s, p}(\R)} = \|\jb{\xi}^s\ft{f}(\xi)\|_{L^p(\R)}$
with $\jb{\cdot} = 1+|\cdot|$.
When $s = 0$, we set $\mathcal{F}L^{p} = \mathcal{F}L^{0, p}$.
The homogeneous $\mathcal{F}\dot{L}^{s, p}$ norm is invariant under the dilation scaling 
when $sp = -1$.

In \cite{VV}, Vargas-Vega 
constructed (both local and global-in-time) solutions for initial data with infinite $L^2$-norm
under certain conditions.
This class of initial data, in particular, contains
those satisfying
\begin{equation} \label{decay}
\bigg|\frac{d^j}{d \xi^j} \ft{u}_0(\xi) \bigg|\lesssim \jb{\xi}^{-\al-j}, \quad j = 0, 1,
\text{ for some }\al > \tfrac{1}{6}.
\end{equation}

\noi
We point out that $u_0$ satisfying \eqref{decay}
is in $\mathcal{F}L^p(\R)$ with $p > \frac{1}{\al}$.
Gr\"unrock \cite{GRUN} considered \eqref{NLS1} on $\R$
with initial data in $\mathcal{F}L^p(\R)$
and proved local well-posedness for $p < \infty$
and global well-posedness for $2< p < \frac{5}{2}$.
The method relies on the Fourier restriction method.
For the global-in-time argument, 
he adapted Bourgain's high-low method \cite{BO2},
where he separated a function in terms of the size of its Fourier coefficient
instead of its frequency size as in \cite{BO2}.

On $\T$, Christ \cite{Christ} applied the power series method
to construct local-in-time solutions (without uniqueness) 
for \eqref{WNLS1} in $\mathcal{F}L^p(\T)$ for $p < \infty$.
Gr\"unrock-Herr \cite{GH} proved the same result
(with uniqueness in a suitable $X^{s, b}$ space)
via the fixed point argument.
A subtraction of $2 \fint |u|^2 dx \, u$ in the nonlinearity in \eqref{WNLS1}
is essential for continuous dependence.
In \cite{Christ}, it is also stated (without proof)
that \eqref{WNLS1} is global well-posed  in $\mathcal{F}L^p$
for sufficiently small (smooth) initial data.

\begin{table} [t]
\begin{center}
 \begin{tabular}{|c|c|c|c|}     \hline
\vphantom{$\Big|$}  & NLS  on $\R$   
& WNLS  on $\T$  &  NLS  on $\T$ \\ \hline 
\vphantom{$\Big|$} 
GWP in $L^2$ & \cite{Tsutsumi} & \cite{BO1} & \cite{BO1} \\ \hline
\vphantom{$\Big|$} 
Ill-posedness below $L^2$ & \cite{KPV5}, \cite{CCT} & \cite{BGT} & \cite{BGT}, \cite{MOLI}  \\ \hline
\vphantom{$\Big|$} 
Well-posedness in $\mathcal{F}L^p$, $p <\infty$ & \cite{GRUN} 
(GWP for $p \in (2, \frac{5}{2}))$ & \cite{Christ}, \cite{GH} & False \cite{Christ}
 \\ \hline
\vphantom{$\Big|$} 
A priori bound for $s \geq -\frac{1}{6}$ & \cite{KochT}, (\cite{CCT3} for $s >-\frac{1}{12}$) & \cite{CHT} 
& Not known  \\ \hline
\vphantom{$\Big|$} 
Weak continuity in $L^2$ & \cite{GM} &  Theorem \ref{THM:1}  & False \cite{MOLI} \\ \hline

\end{tabular}

\vspace{2mm}

\end{center}
\caption{
Corresponding results
for NLS  on $\R$   
and WNLS  on $\T$ (and NLS on $\T$.)}
\label{TAB1}
\end{table}

\subsection{ A priori bound:}
%$\bullet$ {\bf A priori bound:}
Koch-Tataru \cite{KochT} established an a priori bound
on (smooth) solutions for \eqref{NLS1} in $H^s(\R)$ for $s \geq - \frac{1}{6}$
in the form:
given any $M > 0$, there exist $T, C >0$ such that
for any initial $u_0 \in L^2$ with $\|u_0\|_{H^s} \leq M$, 
we have $\sup_{t \in [0, T]} \|u(t)\|_{H^s} \leq C \|u_0\|_{H^s}$,
where $u$ is a solution of \eqref{NLS1} with initial condition $u_0$.
See Christ-Colliander-Tao \cite{CCT3} for a related result.
This result  yields
the existence on weak solutions (without uniqueness).
In the periodic setting, 
Christ-Holmer-Tataru \cite{CHT}
proved the same result for \eqref{WNLS1} when $s \geq -\frac{1}{6}$.
In \cite{KochT}, relating  mKdV and NLS through  modulated wave 
train solutions,  Koch-Tataru
indicate how the regularity $s = -\frac{1}{6}$ 
 arises by associating  mKdV with initial data in $L^2$  to \eqref{NLS1} with
initial data in $H^{-\frac{1}{6}}$ .

\subsection{Weak continuity in $L^2$:}
The Galilean invariance for \eqref{NLS1} yields the critical regularity $s_c^\infty = 0$. 
i.e. the solution map is not uniformly continuous in $H^s$ for $s < s_c^\infty = 0$.
However, it does not imply  that the solution map is not continuous 
in $H^s$ for $s < 0$ (at least on $\R$.)
Heuristically speaking, given $s_0 \in \R$,
one can consider the weak continuity of the solution map in $H^{s_0}$
as an intermediate step between establishing the continuity in 
(the strong topology of) $H^{s_0}$ and 
 proving the continuity in $H^s$ for $s < s_0$.
For example, recall that if $f_n$ converges weakly in $H^{s_0}$, then it converges strongly in $H^s$
for $s < s_0$ (at least in bounded domains.)
Indeed if there is sufficient regularity for the solution map in $H^s$ for some $s < s_0$,
then its weak continuity in $H^{s_0}$ 
can be treated by the approach used in the works 
of Martel-Merle \cite{MM1, MM4}
and Kenig-Martel \cite{KM} related to the asymptotic stability of 
solitary waves. In these works, weak continuity of the flow map plays
a  central role  in the study of  the linearized operator around
the solitary wave and in rigidity theorems.
See Cui-Kenig \cite{CK1} for a nice discussion on this issue.

There are several recent results in this direction.
On $\R$, Goubet-Molinet \cite{GM} proved the weak continuity of the solution map for \eqref{NLS1}
in $L^2(\R)$.
 Cui-Kenig \cite{CK1, CK2} proved the weak continuity 
in the $s_c^\infty$-critical Sobolev spaces
for other dispersive PDEs.
However, on $\T$, Molinet \cite {MOLI} showed that the solution map for \eqref{NLS1}
is not continuous
from $L^2(\T)$ endowed with weak topology
to the space of distributions $(C^\infty(\T))^*$.
This, in particular, implies that 
the solution map for \eqref{NLS1} is not weakly continuous in $L^2(\T)$.

When  considering the Wick ordered cubic NLS \eqref{WNLS1}, 
we remove one of the resonant interactions.
Indeed, we have the following result on the weak continuity of the solution map for \eqref{WNLS1}. 
\begin{theorem}[Weak continuity of WNLS on $L^2(\T)$] \label{THM:1}
Suppose that $u_{0, n}$ converges weakly to $u_0$ in $L^2(\T)$.
Let $u_n$ and $u$ denote the unique global solutions of \eqref{WNLS1}
with initial data $u_{0, n}$ and $u_0$, respectively.
Then, given $T>0$, we have the following.

\begin{itemize}

\item[(a)]
$u_n$ converges weakly to $u$ in $L^4_{T, x} := L^4([-T, T];L^4(\T))$.

\item[(b)]
For any $|t| \leq T$, $u_n(t)$ converges weakly to $u(t)$ in $L^2(\T)$.
Moreover, this weak convergence is uniform for $|t| \leq T$.
i.e. for any $\phi \in L^2(\T)$,
\[ \lim_{n \to \infty} \sup_{|t|\leq T} |\jb{ u_n(t) -u(t), \phi}_{L^2}| = 0.\]
\end{itemize}
\end{theorem}

\noi
We do not expect the weak continuity in the Strichartz space, 
i.e. in $L^6_{T, x}$ (with $|t| \leq T$.)
This is due to the failure of the $L^6_{x, t}$ Strichartz estimate in the periodic setting  \cite{BO1}.
Although the proof of Theorem \ref{THM:1} is essentially contained in \cite{MOLI},
we present it in Section 4 for the completeness of our  presentation.

\section{Wick ordering}

\subsection{Gaussian measures and Hermite polynomials}
In this subsection, we briefly go over the basic relation between Gaussian measures and Hermite polynomials.
For the following discussion,
we refer to  the works of Kuo \cite{KUO2}, Ledoux-Talagrand \cite{Ledoux},
and Janson \cite{Janson}. A nice summary is given by Tzvetkov in
\cite[Section 3]{TZ3} for the hypercontractivity of the Ornstein-Uhlenbeck semigroup related to
products of Gaussian random variables. 

Let $\nu$ be the Gaussian measure with mean 0 and variance $\s$,
and  $H_n(x;\s)$ be the Hermite polynomial of degree $n$ with parameter $\s$. 
They are defined by
\[ e^{tx -\frac{1}{2}\s t^2} = \sum_{n = 0}^\infty \frac{H_n(x;\s)}{n!} t^n.\]

\noi The first three Hermite polynomials are: 
$H_0(x;\s) = 1$, $H_1(x;\s) = x$, and $H_2(x;\s) = x^2-\s$.
It is well known that 
every function $f \in L^2(\nu)$ has a unique series expansion
\[ f(x) = \sum_{n = 0}^\infty a_n \frac{H_n(x;\s)}{\sqrt{n!\s^n}},\]

\noi
where $a_n = (n!\s^n)^{-\frac{1}{2}} \int_{-\infty}^\infty f(x) H_n(x;\s) d\nu(x)$, $n\geq 0$. 
Moreover, we have $\|f\|_{L^2(\nu)}^2 = \sum_{n = 0}^\infty a_n^2.$
In the following, we set $H_n(x) := H_n(x;1)$.

Now, consider the Hilbert space $ L^2(\R^d, \mu_d)$ with $d\mu_d
= (2\pi)^{-\frac{d}{2}} \exp(-{|x|^2}/{2})dx$, $x = (x_1, \dots,
x_d)\in \R^d$. We define a {\it homogeneous Wiener chaos of order
$n$} to be an element of the form $\prod_{j = 1}^d H_{n_j}(x_j)$, $n
= n_1 + \cdots + n_d$. 
Denote  by $\mathcal{K}_n $ the collection of the homogeneous
chaoses of order $n$. Given a homogeneous polynomial $P_n(x) =
P_n(x_1, \dots, x_d)$ of degree $n$, we define {\it the Wick ordered
monomial $:\!P_n(x)\!\!:$} to be its projection onto
$\mathcal{K}_n$. In particular, we have $:x_j^n: = H_{n}(x_j)$ and
$:\prod_{j = 1}^d x_j^{n_j}: = \prod_{j = 1}^d H_{n_j}(x_j)$ with $n
= n_1 + \cdots + n_d$.

In the following, we discuss the key estimate
for the well-posedness results of the Wick ordered cubic NLS of \cite{BO7, CO1}.
Consider the Hartree-Fock operator $L = \Dl -
x \cdot \nabla$, which is the generator for the Ornstein-Uhlenbeck
semigroup. Then, by the hypercontractivity of the Ornstein-Uhlenbeck
semigroup $U(t) = e^{Lt}$, we have the following proposition.

\begin{proposition} 
Fix $q \geq 2$. For every $f \in  L^2(\R^d, \mu_d)$ and $t \geq
\frac{1}{2}\log(q-1)$, we have
\begin{equation}\label{hyp1}
\|U(t) f \|_{L^q(\R^d, \mu_d)}\leq \|f\|_{L^2(\R^d, \mu_d)}.
\end{equation}
\end{proposition}

\noi It is known that the eigenfunction of $L$ with eigenvalue $-n$
is precisely the homogeneous Wiener chaos of order $n$. Thus, we
have the following dimension-independent estimate.

\begin{proposition} \label{PROP:hyp}
Let $F(x)$ be a linear combination of homogeneous chaoses of order
$n$. Then, for $q \geq 2$, we have
\begin{equation} \label{hyp2}
\| F(x) \|_{L^q(\R^d, \mu_d)}\leq (q-1)^\frac{n}{2}
\|F(x)\|_{L^2(\R^d, \mu_d)}.
\end{equation}
\end{proposition}

\noi The proof is basically the same as in \cite[Propositions
3.3--3.5]{TZ3}. We only have to note that $F(x)$ is an eigenfunction
of $U(t) $ with eigenvalue $e^{-nt}$. The estimate  \eqref{hyp2}
follows from \eqref{hyp1} by evaluating \eqref{hyp1} at time $t =
\frac{1}{2} \log (q - 1)$.
In \cite{BO7, CO1, TZ3}, Proposition \ref{PROP:hyp} was used in a crucial manner
to estimate random elements in the nonlinearity
after dyadic decompositions.

In order to motivate $:\!|u|^4\!:$, the Wick ordered $|u|^4$,  for a complex-valued function $u$,
we  consider the Wick ordering on complex Gaussian random variables. 
Let $g$ denote a standard complex-valued Gaussian random
variable. Then, $g$ can be written as $g = x+ iy$, where $x$ and $y$
are independent standard real-valued Gaussian random variables. Note
that the variance of $g$ is $\text{Var}(g)  = 2$. 

Next, we
investigate the Wick ordering on $|g|^{2n}$ for $n\in \mathbb{N}$,
that is, the projection of $|g|^{2n}$ onto $\mathcal{K}_{2n}$.
When $n = 1$, $|g|^2 = x^2 + y^2$ is Wick-ordered into
\[ :|g|^2: = (x^2 - 1) + (y^2-1) = |g|^2 - \text{Var} (g).\]

\noi
When $n = 2$, $|g|^4 = (x^2+ y^2)^2 = x^4 + 2 x^2 y^2 + y^4$ is
Wick-ordered into
\begin{align}
:|g|^4: & = (x^4 -6 x^2 + 3) + 2(x^2 - 1)(y^2 - 1)
+ (y^4 -6 y^2 + 3)\notag \\
& = x^4 + 2 x^2 y^2 + y^4
- 8 (x^2 + y^2) + 8  \label{g4} \\
& = |g|^4 - 4 \text{Var}(g) |g|^2 + 2 \text{Var}(g)^2,\notag
\end{align}

\noi where we used $H_4(x) = x^4 - 6x^2 + 3$.
In general, we have $:|g|^{2n}\!: \, \in \mathcal{K}_{2n}$.
Moreover, we have
\begin{equation}\label{Wick1}
:|g|^{2n}\!: \ = |g|^{2n} + \sum_{j = 0}^{n-1} a_j |g|^{2j}
 = |g|^{2n} + \sum_{j = 0}^{n-1} b_j :|g|^{2j}:  .
\end{equation}

\noi This follows from the fact that $|g|^{2n}$, as a polynomial in
$x$ and $y$ only with even powers, is orthogonal to any homogeneous
chaos of odd order, and it is radial, i.e., it depends only on $|g|^2
= x^2 + y^2$. Note that $:|g|^{2n}\!:$ can also be obtained from the
Gram-Schmidt process applied to $|g|^{2k}$, $ k = 0, \dots, n$ with
$\mu_2 = (2\pi)^{-1} \exp (-(x^2 + y^2)/2) dx dy$.

\subsection{Wick ordered cubic NLS}
In \cite{BO7}, Bourgain considered the issue of the invariant Gibbs measure
for \eqref{NLS1} on $\T^2$ in the defocusing case.
In this subsection, we present his argument to derive \eqref{WNLS1} on $\T^2$.
First, consider the finite dimensional approximation to \eqref{NLS1}:
\begin{equation} \label{NLS2}
\begin{cases}
i u^N_t - \Dl u^N + \mathbb{P}_N (|u^N|^2 u^N) = 0\\
u\big|_{t = 0} = \mathbb{P}_N u_0,
\end{cases}
\quad (x, t) \in \T^2\times \R,
\end{equation}

\noi
where  $u^N=\mathbb{P}_N u$ and 
 $\mathbb{P}_N$ is the Dirichlet projection onto the frequencies $|n|\leq N$.
This is a Hamiltonian equation with Hamiltonian
$H(u^N)$, where $H$ is as in \eqref{conserved} with the $+$ sign. 
On $\T^2$, the Gaussian part $d \rho = Z^{-1} \exp \big(-\frac{1}{2} \int |\nabla u|^2 dx \big) \prod_{x\in\T^2} d u(x)$
of the Gibbs measure is supported on $\bigcap_{s <0} H^s (\T^2) \setminus L^2(\T^2)$.
However, the nonlinear part $\int |\mathbb{P}_N u|^4 dx$
of the Hamiltonian diverges to $\infty$ as $N\to \infty$
almost surely on the support of the Wiener measure $\rho$.
Hence, we need to renormalize the nonlinearity.

A typical element  in the support of the Wiener measure $\rho$ is given by
\begin{equation}\label{IV3}
u (x;\omega) = \sum_{n\in\Z^2} \frac{g_n(\omega)}{\sqrt{1+|n|^{2}}}e^{in\cdot x}, 
\end{equation}

\noi
where $\{g_n\}_{n\in\Z}$ is a family of independent standard 
complex-valued Gaussian random variables.\footnote{The expression \eqref{IV3}
is a representation of elements in the support of 
$d \wt{\rho} = \wt{Z}^{-1} \exp \big(-\frac{1}{2} \int |u|^2 -\frac{1}{2} \int |\nabla u|^2  \big) \prod_{x\in\T^2} d u(x)$
due to the problems at the zero Fourier mode for $\rho$.
However, we do not worry about this issue.}
For simplicity, assume that $\text{Var}(g_n) = 1$.
For $u$ of the form \eqref{IV3}, define $a_N$ by 
\[a_N = \mathbb{E} \Big[ \int \hspace{-11pt} - |u^N|^2 dx \Big] = \sum_{|n|\leq N} \frac{1}{1+|n|^2}.\]

\noi
We have that  $a_N \sim \log N$ for large $N$.
We define the Wick ordered truncated Hamiltonian $H_N$ by
\begin{align} \label{u4}
H_N (u^N) &= \frac{1}{2} \int_{\T^2} |\nabla u^N|^2 dx +\frac{1}{4}\int_{\T^2} :| u^N|^4\!: dx\\
 &= \frac{1}{2} \int_{\T^2} |\nabla u^N|^2 dx +\frac{1}{4}\int_{\T^2} | u^N|^4 -4 a_N |u^N|^2 + 2a_N^2 \ dx. \notag
\end{align}

\noi
(Compare \eqref{u4} with \eqref{g4}.)
From  \eqref{u4}, we obtain an Hamiltonian equation that is 
the Wick ordered version of \eqref{NLS2}:
\begin{equation} \label{WNLS2}
i u^N_t - \Dl u^N + \mathbb{P}_N (|u^N|^2 u^N) - 2 a_N u^N= 0.
\end{equation}

\noi
Let $c_N = \fint|u^N|^2 - a_N$, we see that 
$c_\infty (\omega) = \lim_{N\to \infty} c_N (\omega) < \infty$ almost surely.
Under the change of variables  $v^N = e^{-2i c_N t} u^N$,  \eqref{WNLS2} becomes
\begin{equation} \label{WNLS3}
\textstyle i v^N_t - \Dl v^N + \mathbb{P}_N \big(|v^N|^2 -2\fint |v^N|^2\big) v^N= 0.
\end{equation}

\noi
Finally, letting $N\to \infty$, we  obtain the Wick order NLS.
\begin{equation} \label{WNLS4}
\textstyle  i v_t - \Dl v + (|v|^2 -2\fint |v|^2\big) v= 0.
\end{equation}

On $\T$, one can repeat the same argument.
Note the following issue.
On the one hand, the assumption that $u(t)$ is of the form \eqref{IV}
is natural for $\al \in \mathbb{N} \cup \{0\}$
in view of the conservation laws.
On the other hand, 
$c_N = \fint|u^N|^2 - \mathbb{E}\big[\fint|u^N|^2\big] <\infty$
for $\al > \frac{1}{4}$.
i.e. $\al = 1$ is the smallest integer value of such $\al$.
In this case, there is no need for the Wick ordered NLS \eqref{WNLS1}
since $u \in H^{\frac{1}{2}-}$ a.s. for $\al = 1$. 

\section{Weak continuity of the Wick ordered cubic NLS in $L^2(\T)$}

In this section, we present the proof of Theorem \ref{THM:1}.
First, write \eqref{WNLS1}   as an integral equation:
\begin{equation}\label{WWNLS2}
u(t) = S(t) u_0 \pm i \int_0^t S(t - t') \mathcal{N}(u)(t') dt'
\end{equation}

\noi
where $\mathcal{N}(u)= (|u|^2 - 2 \fint |u|^2) u$
and $S(t) = e^{-i \dx^2 t}$.
Define $\mathcal{N}_1(u_1, u_2, u_3)$ and $\mathcal{N}_2(u_1, u_2, u_3)$
by
\begin{align*}
\mathcal{N}_1(u_1, u_2, u_3)
& = \sum_{\substack{n = n_1 - n_2 + n_3\\
n_2 \ne n_1, n_3}}
\ft{u}_1(n_1)\cj{\ft{u}}_2(n_2)\ft{u}_3(n_3)e^{inx},\\
\mathcal{N}_2(u_1, u_2, u_3)
& = - \sum_{n}
\ft{u}_1(n)\cj{\ft{u}}_2(n)\ft{u}_3(n)e^{inx}.
\end{align*}

\noi
Moreover, let $\mathcal{N}_j(u) := \mathcal{N}_j(u, u, u)$.
Then, we have
$\mathcal{N}(u)= \mathcal{N}_1(u)+\mathcal{N}_2(u).$

In \cite{BO1}, Bourgain established global well-posedness of \eqref{NLS1} (and \eqref{WNLS1})
by introducing a new weighted space-time Sobolev space $X^{s, b}$
whose norm is given by
\[ \|u\|_{X^{s, b}(\T \times\R)} = \|\jb{n}^s \jb{\tau - n^2}^b \ft{u}(n, \tau) \|_{L^2_{n, \tau}(\Z\times \R)}\] 

\noi
where $\jb{\,\cdot\,} = 1 + |\cdot|$.
Define the local-in-time version $X^{s, b}_\dl$ on $[-\dl, \dl]$
by
\[ \|u\|_{X^{s, b}_\dl} = \inf \big\{ \|\wt{u}\|_{X^{s, b}}\, ; \wt{u}|_{[-\dl, \dl]} = u \big\}.\]

\noi
In the following, we list
the estimates needed for local well-posedness of \eqref{WNLS1}.
Let $\eta(t)$ be a smooth cutoff function such that $\eta = 1$ on $[-1, 1]$
and $\eta = 0$ on $[-2, 2]$.
\begin{itemize}
\item Homogeneous linear estimate: for $s, b \in \R$, we have
\begin{equation} \label{lin1}
\| \eta(t) S(t) f \|_{X^{s, b}} \leq C_1 \| f\|_{H^s}. 
\end{equation}

\item Nonhomogeneous linear estimate: for $s \in \R$ and $b > \tfrac{1}{2}$, we have
\begin{equation} \label{lin2}
\bigg\|\eta(t) \int_0^t S(t - t') F(t') dt' \bigg\|_{X_\dl^{s, b}}
\lesssim C(\dl) \|F\|_{X_\dl^{s, b-1}}.
\end{equation}

\item Periodic $L^4$ Strichartz estimate:
Zygmund \cite{ZYG} proved
\begin{equation} \label{L41}
\| S(t) f \|_{L^4_{x, t} (\T\times [-1, 1])} \lesssim \|f\|_{L^2},
\end{equation}

\noi
which was improved by Bourgain \cite{BO1}: 
\begin{equation} \label{L42}
\| u \|_{L^4_{x, t} (\T\times [-1, 1])} \lesssim \|u\|_{X^{0, \frac{3}{8}}}.
\end{equation}

\end{itemize}

\noi
These estimates allow us to prove local well-posedness of  \eqref{WNLS1}
via the fixed point argument
such that a solution $u$ exists on the time interval $[-\dl, \dl]$
with $\dl = \dl(\|u_0\|_{L^2})$.
Moreover, we have $\|u\|_{X_\dl^{0,\frac{1}{2}+}} \lesssim  \|u_0\|_{L^2}$.
Such local-in-time solutions can be extended globally in time thanks to the $L^2$ conservation.

\medskip

Now, fix $u_0 \in L^2(\T)$, and  let $u_{0, n}$ converges weakly to $u_0$ in $L^2(\T)$.
Denote by $u_n$ and $u$
the unique global solutions of \eqref{WNLS1}
with initial data $u_{0, n}$ and $u_0$. 
By the uniform boundedness principle, 
we have $\|u_{0, n}\|_{L^2}, \, \|u_0\|_{L^2} \leq C$ for some $C>0$.
Hence, the local well-posedness guarantees 
the existence of the solutions $u_n$, $u$ on the time interval $[-\dl, \dl]$
with $\dl = \dl (C)$, uniformly in $n$.
In the following, we assume $\dl = 1$.
i.e. we assume that all the estimates hold on $[-1, 1]$.
(Otherwise we can replace $[-1, 1]$ by $[-\dl, \dl]$ for some $\dl >0$
and iterate the argument in view of the $L^2$ conservation.)

\subsection{Proof of Theorem \ref{THM:1} (a)}

First, we show that $u_n$ converges to $u$ as space-time distributions.

\noi
$\bullet$ {\bf Linear part:}
Since $u_{0, n} \wto u_0$ in $L^2(\T)$, 
we have $\| u_{0, n} - u_0\|_{H^{-\eps}(\T)} \to 0$ for any $\eps > 0$.
Let $\phi \in C_c^\infty(\T\times \R)$ be a test function. 
Then, by H\"older inequality and \eqref{lin1}, we have
\begin{align*}
\iint  \eta(t) S(t)(u_{0, n} - u_0) \phi(x, t) dx dt
& \leq \|\eta(t)  S(t) (u_{0, n} - u_0)\|_{X^{-\eps, \frac{1}{2}+}}
\|\phi\|_{X^{\eps, -\frac{1}{2}-}} \\
& \lesssim C_\phi \|u_{0, n} - u_0\|_{H^{-\eps}} \to 0.
\end{align*}

\noi
Hence, $\eta(t)S(t)u_{0, n}$ converges to $\eta(t)S(t)u_{0}$ as space-time distributions.

\medskip

\noi
$\bullet$ {\bf Nonlinear part:}
Let $\mathcal{M}(u)$ denote the Duhamel term.
i.e.
\[\mathcal{M}(u)(t):= \pm i \int_0^t S(t - t') \mathcal{N}(u)(t') dt'.\]

\noi
Similarly, define
$\mathcal{M}_j(u_1, u_2, u_3)$
by 
\[\mathcal{M}_j(u_1, u_2, u_3)(t):= \pm i \int_0^t S(t - t') \mathcal{N}_j(u_1, u_2, u_3)(t') dt'\]

\noi
for $j = 1, 2$. 
Also, let
$\mathcal{M}_j(u):= \mathcal{M}_j(u, u, u)$.

From  the local theory, we have
$\|u_n\|_{X_1^{0, \frac{1}{2}+} }\lesssim \|u_{0, n}\|_{L^2} \leq C$ for all $n$.
Thus, there exists a subsequence 
$u_{n_k}$ converging weakly to some $v$ in $X_1^{0, \frac{1}{2}+}$.
It  follows from \cite[Lemmata 2.2 and 2.3]{MOLI}
that 
$\mathcal{N}_j$, $j = 1, 2$, 
is weakly continuous from $X_1^{0, \frac{1}{2}+}$ into $X_1^{0, -\frac{7}{16}}$.
Hence, $\mathcal{N}_j (u_k) \wto \mathcal{N}_j (v)$ in $X_1^{0, -\frac{7}{16}}$.

Recall the following.
Given  Banach spaces $X$ and $Y$ with a continuous linear operator $T: X\to Y$, 
we have $T^*: Y^* \to X^*$.
If $f_n \wto f$ in $X$, then we have, for $\phi \in Y^*$, 
$ \jb{T(f_n-f), \phi} = \jb{ f_n - f, T^* \phi} \to 0$
since $T^* \phi \in X^*$.
Hence, $Tf_n \wto Tf$ in $Y$.

It follows from \eqref{lin2} that  the map: $F \longmapsto \int_0^t S(t - t') F(t') dt'$
is linear and continuous from $X_1^{0, -\frac{7}{16}}$
into $X_1^{0, \frac{1}{2}+}$.
Hence,  $\mathcal{M}(u_{n_k}) \wto  \mathcal{M}(v)$ in $X_1^{0, \frac{1}{2}+}$.
In particular, 
$\mathcal{M}(u_{n_k})$  converges to $\mathcal{M}(v)$
as space-time distributions.

Since $u_{n_k}$ is a solution to \eqref{WNLS1} with  initial data $u_{0, n_k}$, we have
\[ u_{n_k} = \eta \, S(t) u_{0, n_k} + \eta \,  \mathcal{M}(u_{n_k}).\]

\noi
By taking the limits of both sides, we obtain
\begin{equation*} 
 v = \eta \,  S(t)u_0 + \eta \,  \mathcal{M}(v),
\end{equation*}

\noi
where the equality holds in the sense of space-time distributions.
From the uniqueness of solutions to \eqref{WNLS1} in $X^{0, \frac{1}{2}+}_1$, 
we have $v = u$ in $X_1^{0, \frac{1}{2}+}$.

In fact, we can show that uniqueness of solutions to \eqref{WNLS1}
holds in $L^4_{x, t}(\T \times [-1, 1])$ with little effort.
For simplicity, we replace $\mathcal{N}(u)$ in \eqref{WWNLS2} by $|u|^2 u$.
Then, by \eqref{L41} and \eqref{L42}, we have
\begin{align*}
\|\eta(t) u\|_{L^4_{x, t}} 
& \leq \|\eta(t) S(t) u_0 \|_{L^4_{x, t}}
+ \bigg\|\eta(t) \int_0^t S(t - t') |\eta u (t')|^2 \eta u(t')dt' \bigg\|_{L^4_{x, t}}\\
& \lesssim \|u_0\|_{L^2_x} 
+ \bigg\|\eta(t) \int_0^t S(t - t') |\eta u(t')|^2 \eta u(t')dt' \bigg\|_{X^{0, \frac{3}{8}}}.
\end{align*}

\noi
Moreover, 
we can use \eqref{lin2}, duality,  $L^4_{x, t}L^4_{x, t}L^4_{x, t}L^4_{x, t}$-H\"older inequality,
and \eqref{L42}
to estimate the second term by
\begin{align*}
& \lesssim  \| |\eta u|^2 \eta u\|_{X^{0, -\frac{3}{8}}}
=  \sup_{\|v\|_{X^{0, \frac{3}{8}}}= 1} \iint v |\eta u|^2 (\eta u) dx dt
\leq  \sup_{\|v\|_{X^{0, \frac{3}{8}}}= 1} \|v\|_{L^4_{x, t}}\|\eta u\|^3_{L^4_{x, t}}
\leq \|\eta u\|^3_{L^4_{x, t}}.
\end{align*}

\noi
This shows that $u$ is indeed a unique solution in $L^4_{x, t}(\T\times [-1, 1])$

It follows from $(L^4_{x, t}(\T \times [-1, 1]))^* \subset (X^{0, \frac{1}{2}+}_1)^*$
that weak convergence in $X^{0, \frac{1}{2}+}_1$ implies weak convergence 
in $L^4_{x, t}(\T \times [-1, 1])$.
Hence, the subsequence $u_{n_k}$ converges weakly to $u$ in 
$X^{0, \frac{1}{2}+}_1$ and $L^4_{x, t}(\T \times [-1, 1])$.
The argument above also shows that $u$ is the only weak limit point of $u_n$
in $X^{0, \frac{1}{2}+}_1$ and $L^4_{x, t}(\T \times [-1, 1])$.
Then,  it follows from the boundedness of $u_n$ in $X^{0, \frac{1}{2}+}_1$ and $L^4_{x, t}(\T \times [-1, 1])$
that the whole sequence $u_n$ converges weakly to $u$.
Indeed, suppose that the whole sequence $u_n$ does not converge weakly to $u$.
Then, there exists $\phi \in (X^{0, \frac{1}{2}+}_1)^*$ 
such that $\jb{ u_n, \phi} \nrightarrow \jb{u, \phi}$.
This, in turn, implies that
there exists $\eps > 0$ such that for any $N \in \mathbb{N}$,
there exists $n \geq N$ such that 
$|\jb{ u_n - u , \phi}| > \eps$.
Given $\eps > 0$, we can construct a subsequence $u_{n_k}$
with $|\jb{ u_{n_k} - u , \phi}| > \eps$ for each $k$.
However, by repeating the previous argument (from the uniform boundedness of $u_{n_k}$
in $X^{0, \frac{1}{2}+}_1$), 
$u_{n_k}$ has a sub-subsequence converging to $u$, which is a contradiction.
This establishes Part (a) of Theorem \ref{THM:1}
on $[-1, 1]$.

\subsection{Proof of Theorem \ref{THM:1} (b)}
Recall the following embedding. 
For $b >\frac{1}{2}$, we have
\begin{equation} \label{lin3}
\|u\|_{L^\infty([-1, 1]; H^s)} \leq C_2 \|u\|_{X^{s, b}_1}.
\end{equation}

\noi
Fix $\phi \in L^2(\T)$ in the following.

\medskip

\noi
$\bullet$ {\bf Linear part:}
Given $\eps> 0$, choose $\psi \in H^1(\T)$ 
such that $\| \phi - \psi \|_{L^2} < \frac{\eps}{2KC_1C_2}$,
where $K = \sup_n \| u_{0, n} - u_0\|_{L^2} <\infty$
and $C_1$, $C_2$ are as in \eqref{lin1}, \eqref{lin3}.
Then, by \eqref{lin1}
and \eqref{lin3}, we have
\begin{align*}
\sup_{|t| \leq 1} | \jb{ S(t)(u_{0, n} - u_0, \phi}_{L^2}|
& \leq \sup_{|t| \leq 1} | \jb{ S(t)(u_{0, n} - u_0, \psi}_{L^2}|
+ \sup_{|t| \leq 1} | \jb{ S(t)(u_{0, n} - u_0, \phi - \psi}_{L^2}|\\
& \leq \| S(t)(u_{0, n} - u_0)\|_{L^\infty([-1, 1]; H^{-1})} \|\psi\|_{H^1}\\
& \hphantom{XXXXXXX} +  \| S(t)(u_{0, n} - u_0)\|_{L^\infty([-1, 1]; L^2)} \|\phi - \psi\|_{L^2}\\
& \leq C_\psi\| S(t)(u_{0, n} - u_0)\|_{X_1^{-1, \frac{1}{2}+}}
+ \tfrac{\eps}{2KC_1} \| S(t)(u_{0, n} - u_0)\|_{X_1^{0, \frac{1}{2}+}} \\
& \leq C \| u_{0, n} - u_0\|_{H^{-1}}
+ \tfrac{\eps}{2K}\| u_{0, n} - u_0\|_{L^2}.
\end{align*}

\noi
Hence, there exists $N_1$ such that for $n \geq N_1$,
\[\sup_{|t| \leq 1} | \jb{ S(t)(u_{0, n} - u_0, \phi}_{L^2}| <\eps\]

\noi
since $u_{0, n}$ converges strongly $u_n$ in $H^{-1}$.

\medskip

\noi
$\bullet$ {\bf Nonlinear part:}
Since $u_n \wto u$ in $X_1^{0, \frac{1}{2}+}$, 
we see that 
$\mathcal{N}(u_n)$ converges strongly to $\mathcal{N}(u)$
in $X_1^{-1, - \frac{7}{16}}$.
See the proof of Lemmata 2.2 and 2.3 in \cite{MOLI}.
Then, it follows from \eqref{lin2}
that  $\mathcal{M}(u_n)$ converges strongly to  $\mathcal{M}(u)$
in $X_1^{-1, \frac{1}{2}+}$.

Given $\eps> 0$, choose $\psi \in H^1(\T)$ such that 
$\| \phi - \psi \|_{L^2} < \frac{\eps}{2KC_2}$,
where $K = \sup_n \|\mathcal{M}(u_n) - \mathcal{M}(u_n)\|_{X_1^{0, \frac{1}{2}+}} <\infty$
and $C_2$ is as in \eqref{lin3}.
Then, by \eqref{lin3}, we have
\begin{align*}
\sup_{|t| \leq 1} & | \jb{ \mathcal{M}(u_n) - \mathcal{M}(u), \phi}|\\
& \leq \sup_{|t| \leq 1} | \jb{ \mathcal{M}(u_n) - \mathcal{M}(u), \psi}_{L^2}| 
+ \sup_{|t| \leq 1} | \jb{ \mathcal{M}(u_n) - \mathcal{M}(u), \phi - \psi}_{L^2}|\\
& \leq  \|\mathcal{M}(u_n) - \mathcal{M}(u)\|_{L^\infty([-1, 1]; H^{-1})} \|\psi\|_{H^1}\\
&  \hphantom{XXXXX} +  \|\mathcal{M}(u_n) - \mathcal{M}(u)\|_{L^\infty([-1, 1]; L^2)} \|\phi - \psi\|_{L^2}\\
& \leq C_\psi \|\mathcal{M}(u_n) - \mathcal{M}(u)\|_{X_1^{-1, \frac{1}{2}+}}
+ \tfrac{\eps}{2K} \|\mathcal{M}(u_n) - \mathcal{M}(u)\|_{X_1^{0, \frac{1}{2}+}}.
\end{align*}

\noi
Hence, there exists $N_2$ such that
for $n \geq N_2$,
\[\sup_{|t| \leq 1}  | \jb{ \mathcal{M}(u_n) - \mathcal{M}(u), \phi}| <\eps.\]

\medskip

\noi
Therefore, we have 
\begin{equation*}
 \lim_{n \to \infty} \sup_{|t|\leq 1} |\jb{ u_n(t) -u(t), \phi}_{L^2}| = 0.
\end{equation*}

\noi
Given $[-T, T]$, we can iterate Part 1 and 2 on each $[j, j+1]$
and obtain Theorem \ref{THM:1}.

\section*{Acknowledgements}
We would like to  thank Jim Colliander for helpful discussions and comments
and Luis Vega for pointing out additional references.

\end{document}